\documentclass{article}
\usepackage{amsmath,amsthm,amsopn,amstext,amscd,amsfonts,amssymb}
\usepackage{natbib}

\def\diag{\mathop{\rm diag}\nolimits}
\def\tr{\mathop{\rm tr}\nolimits}

\def\build#1#2#3{\mathrel{\mathop{#1}\limits^{#2}_{#3}}}

\def\etr{\mathop{\rm etr}\nolimits}

\renewenvironment{abstract}
                 {\vspace{6pt}
                  \begin{center}
                  \begin{minipage}{5in}
                  \centerline{\textbf{Abstract}}
                  \noindent\ignorespaces
                 }
                 {\end{minipage}\end{center}}
\newtheorem{thm}{\textbf{Theorem}}[section]
\newtheorem{cor}{\textbf{Corollary}}[section]
\newtheorem{lem}{\textbf{Lemma}}[section]
\theoremstyle{definition}

\newtheorem{rem}{\textbf{Remark}}[section]

\title{\Large \textbf{Shape theory via polar  decomposition}}
\author{
  \textbf{Jos\'e A. D\'{\i}az-Garc\'{\i}a} \thanks{Corresponding author\newline
    }\\
  {\normalsize Department of Statistics and Computation} \\
  {\normalsize Universidad Aut\'onoma Agraria Antonio Narro}\\
  {\normalsize 25350 Buenavista, Saltillo, Coahuila, Mexico} \\
  {\normalsize E-mail: jadiaz@uaaan.mx} \\[2ex]
  \textbf{Francisco J. Caro-Lopera} \\
  {\normalsize Department of Basic Sciences} \\
  {\normalsize Universidad de Medell\'{\i}n} \\
  {\normalsize Carrera 87 No.30-65, of. 5-103}\\
  {\normalsize Medell\'{\i}n, Colombia}\\
  {\normalsize E-mail: fjcaro@udem.edu.co}\\
}
\date{}
\begin{document}
\maketitle

\begin{abstract}
This work proposes a new model in the context of statistical theory of shape, based
on the polar decomposition. The non isotropic noncentral elliptical shape
distributions via polar decomposition is derived in the context of zonal polynomials,
avoiding the invariant polynomials and the open problems for their computation. The
new polar shape distributions are easily computable and then the inference procedure
can be studied under exact densities. As an example of the technique, a classical
application in Biology is studied under three models, the usual Gaussian and two non
normal Kotz models; the best model is selected by a modified BIC criterion, then a
test for equality in polar shapes is performed.
\end{abstract}

\section{Introduction}

Matrix variate statistical shape analysis has been extensively studied in the last
two decades by a number of approaches: via QR decomposition (\citet{GM93}); via SVD
decompositions (\citet{g:91}, \citet{LK93}, \citet{dgm:97}, \citet{dgr:03}); via
affine transformations (\citet{GM93}, \citet{dgr:03}, \citet{Caro2009}), among many
others methods (\citet{DM98} and the references there in). However, the polar
decomposition has not been included yet in the context of shape theory.

According to the transformation, we say that the shape of an object is all
geometrical information which remains after filtering out translation, scale,
rotation, reflection, uniform share, etc, from an original figure comprised in $N$
landmarks in $K$ dimension. Statistical shape theory study the mean shape of
populations in presence of randomness.

Some of the classical works (\citet{GM93}) assume an isotropic Gaussian model for the
landmark matrix in order to obtain shape densities expanded in known polynomials,
such as zonal polynomials (\citet{JAT64}, \citet{MR1982}); then, generalisations for
matrix variate shape theory under elliptical models appeared via the SVD method
(\citet{dgr:03}) and via the affine technique (\citet{Caro2009}). However, in order
to obtain zonal polynomials a partial non isotropy was assumed, otherwise,
considering a full non isotropy, the densities are expanded in terms of invariant
polynomials (\citet{D80}), which are non available for large degrees.

Now, the isotropic assumption, say $\boldsymbol{\Theta} = \mathbf{I}_{K}$ for an
elliptical shape model of the form
$$
   \mathbf{X} \sim \mathcal{E}_{N \times K} (\boldsymbol{\mu}_{{}_{\mathbf{X}}},
   \boldsymbol{\Sigma}_{{}_{\mathbf{X}}},\boldsymbol{\Theta}, h),
$$
restricts substantially the correlations of the landmarks in the figure and it is non
appropriate for applications.  So, we expect the non isotropic model, with any
positive definite matrix $\boldsymbol{\Theta}$, as the best model for considering all
the possible correlations among the anatomical (geometrical o mathematical) points.

This work solves that problem and sets the non isotropic noncentral elliptical shape
distributions via polar decomposition in the context of zonal polynomials, avoiding
the invariant polynomials and the open problems for their computation. The new shape
distributions are easily computable and then the inference procedure can be studied
under exact densities.

In section \ref{sec:polarsizeandshape}, the so termed polar shape coordinates are
introduced and the main mathematical tools are studied in order to obtained the polar
size and shape density. Then the polar shape density is derived in section
\ref{sec:polarshape}. Section \ref{sec:central} studies the central case and the
corresponding invariance under the family of elliptical distributions. Finally,
section \ref{sec:particularmodelsexample} gives explicit densities and performs
inference with three models, the Gaussian and two non normal Kotz models.

\section{Polar size-and-shape
distribution}\label{sec:polarsizeandshape}

Consider a full non isotropy (non singular) elliptical model
$$
   \mathbf{X} \sim \mathcal{E}_{N \times K} (\boldsymbol{\mu}_{{}_{\mathbf{X}}},
   \boldsymbol{\Sigma}_{{}_{\mathbf{X}}},\boldsymbol{\Theta}, h),
$$
with  generator function $h(\cdot)$.

In order to avoid the referred problem of invariant polynomials consider the
following procedure: Let
$$
   \mathbf{X} \sim {\mathcal E}_{N \times K} (\boldsymbol{\mu}_{{}_{\mathbf{X}}},
   \boldsymbol{\Sigma}_{{}_{\mathbf{X}}},\boldsymbol{\Theta}, h),
$$
if $\boldsymbol{\Theta}^{1/2}$ is the positive definite square root of the matrix
$\boldsymbol{\Theta}$, i .e. $\boldsymbol{\Theta} = (\boldsymbol{\Theta}^{1/2})^{2}$,
with $\boldsymbol{\Theta}^{1/2}:$ $K \times K$, \citet[p. 11]{gv:93}, and noting that
$$
  \mathbf{X} \boldsymbol{\Theta}^{-1} \mathbf{X}' = \mathbf{X} (\boldsymbol{\Theta}^{-1/2}
  \boldsymbol{\Theta}^{-1/2})^{-1}\mathbf{X}' = \mathbf{X} \boldsymbol{\Theta}^{-1/2} (\mathbf{X}
  \boldsymbol{\Theta}^{-1/2})' = \mathbf{Z}\mathbf{Z}',
$$ where
$$
\mathbf{Z} = \mathbf{X} \boldsymbol{\Theta}^{-1/2},
$$
then
$$
  \mathbf{Z} \sim {\mathcal E}_{N \times K}(\boldsymbol{\mu}_{{}_{\mathbf{Z}}},
  \boldsymbol{\Sigma}_{{}_{\mathbf{X}}}, \mathbf{I}_{K}, h)
$$
with $\boldsymbol{\mu}_{{}_{\mathbf{Z}}} = \boldsymbol{\mu}_{{}_{\mathbf{X}}}
\boldsymbol{\Theta}^{-1/2}$, see \citet[p. 20]{gv:93}.

And we arrive at the classical starting point in shape theory where the original
landmark matrix is replaced by $\mathbf{Z} = \mathbf{X} \boldsymbol{\Theta}^{-1/2}$
(see \citet{GM93}, for example). Then we can proceed as usual, removing from
$\mathbf{Z}$, translation, scale, rotation  in order to obtain the shape of
$\mathbf{Z}$ (or $\mathbf{X}$) via QR, SVD, or  polar decompositions, for example.

In this paper we consider a new system of shape coordinates,  the polar shape
coordinates $\mathbf{u}$ of $\mathbf{X}$ which are constructed as follows:
$$
  \mathbf{L}\mathbf{X}=\mathbf{Y}=\mathbf{R}\mathbf{H}=r\mathbf{W}\mathbf{H} = r
 \mathbf{W}(\mathbf{u})\mathbf{H}.
$$
The matrix $\mathbf{L}$ is as usual an $(N-1)\times N$ Helmert submatrix and we
assume that
$$
  \mathbf{Y}\sim\mathcal{E}_{N-1\times K}(\boldsymbol{\mu}, \boldsymbol{\Sigma}\otimes
  \mathbf{I}_{K},h),\quad \boldsymbol{\mu}=\mathbf{L}\boldsymbol{\mu}_{\mathbf{X}},
  \quad \boldsymbol{\Sigma}=\mathbf{L}\boldsymbol{\Sigma}_{\mathbf{X}}\mathbf{L}'.
$$
Under this approach, $\mathbf{Y}=\mathbf{R}\mathbf{H}$ is the polar decomposition,
where $\mathbf{R}:N-1\times N-1$ is a positive definite matrix and $\mathbf{H}\in
V_{N-1,K}$.

It is important to note that only under $n=N-1$ the polar approach is valid.

Recall that for the singular value decomposition $\mathbf{Y} = \mathbf{P} \mathbf{L}
\mathbf{Q}'$, the polar decomposition of $\mathbf{Y}$ is given by $\mathbf{Y} =
\mathbf{R} \mathbf{H}$, where $\mathbf{R} = \mathbf{P} \mathbf{L }\mathbf{P}'$ and
$\mathbf{H}=\mathbf{P}\mathbf{Q}'$.

So, we start with a known result, see \citet{c:96}:

\begin{lem}\label{lem:polarjacobian}
Let $\mathbf{Y}:N-1\times K$, then there exist $\mathbf{R}:N-1\times N-1$ a positive
definite matrix and $\mathbf{H}\in V_{N-1,K}$ such that $\mathbf{Y} = \mathbf{R}
\mathbf{H}$ and
$$
  (d\mathbf{Y})=\displaystyle\prod_{i<j}^{N-1}(L_{i}+L_{j})(d\mathbf{R})(\mathbf{H}d\mathbf{H}'),
$$
with $\mathbf{L}=\diag (L_{1},\ldots,L_{N-1})$ and
$\mathbf{R}=\mathbf{P}\mathbf{L}\mathbf{P}'$ i.e. $L_{i}=\lambda_{i}(\mathbf{R})$.
\end{lem}

So the main result of the section follows:

\begin{thm}\label{th:Polarsizeandshape}
The polar  size-and-shape density is
\begin{small}
\begin{equation}\label{eq:Polarsizeandshape}
    f_{\mathbf{R}}(\mathbf{R}) = \frac{\pi^{\frac{(N-1)K}{2}}\displaystyle
    \prod_{i<j}^{N-1}(L_{i}+L_{j})}{2^{-N+1}\Gamma_{N-1}\left[\frac{K}{2}\right]
    |\boldsymbol{\Sigma}|^{\frac{K}{2}}} \sum_{t=0}^{\infty}\sum_{\kappa}
    \frac{h^{(2t)}\left[\tr\left(\boldsymbol{\Sigma}^{-1}\mathbf{R}^{2} +
    \boldsymbol{\Omega}\right)\right]}{t!} \frac{C_{\kappa}\left(\boldsymbol{\Omega}
    \boldsymbol{\Sigma}^{-1}\mathbf{R}^{2}\right)}{\left(\frac{1}{2}K\right)_{\kappa}}.
\end{equation}
\end{small}
\end{thm}
\textit{Proof.} The density of $\mathbf{Y}$ with $\boldsymbol{\boldsymbol{\Omega}} =
\boldsymbol{\boldsymbol{\Sigma}}^{-1} \boldsymbol{\mu}
\boldsymbol{\boldsymbol{\Theta}}^{-1} \boldsymbol{\mu}'$ is
$$
  f_{\mathbf{Y}}(\mathbf{Y})=\frac{1}{|\boldsymbol{\Sigma}|^{\frac{K}{2}}}h
  \left[\tr\left(\boldsymbol{\Sigma}^{-1}\mathbf{Y}\mathbf{Y}'+\boldsymbol{\Omega}\right)
  - 2\tr\boldsymbol{\Sigma}^{-1}\mathbf{Y}\boldsymbol{\mu}'\right].
$$
If  the decomposition $\mathbf{Y}=\mathbf{R}\mathbf{H}$ is performed
and the Lemma \ref{lem:polarjacobian} is applied, then the joint
density of $\mathbf{R}$ and $\mathbf{H}$ remains
$$
  f_{\mathbf{R},\mathbf{H}}(\mathbf{R},\mathbf{H})=\frac{\displaystyle
  \prod_{i<j}^{N-1}(L_{i}+L_{j})}{|\boldsymbol{\Sigma}|^{K/2}}h\left[\tr
  \left(\boldsymbol{\Sigma}^{-1}\mathbf{R}^{2}+\boldsymbol{\Omega}\right)-2
  \tr\boldsymbol{\mu}'\boldsymbol{\Sigma}^{-1}\mathbf{R}\mathbf{H}\right].
$$
Assuming that $h(\cdot)$ can be expanded in a convergent power
series, see \citet{fz:90}, i.e.
$$
  h(a+v)=\sum_{t=0}^{\infty}\frac{h^{(t)}(a)}{t!}v^{t},
$$
hence
$$
  f_{\mathbf{R},\mathbf{H}}(\mathbf{R},\mathbf{H})=\frac{\displaystyle
  \prod_{i<j}^{N-1}(L_{i}+L_{j})}{|\boldsymbol{\Sigma}|^{K/2}}\sum_{t=0}^{\infty}
  \frac{h^{(t)}\left[\tr\left(\boldsymbol{\Sigma}^{-1}\mathbf{R}+\boldsymbol{\Omega}
  \right)\right]}{t!}\left[\tr(-2\boldsymbol{\mu}'\boldsymbol{\Sigma}^{-1}\mathbf{R}
  \mathbf{H})\right]^{t}(\mathbf{H}d\mathbf{H}').
$$
So, the marginal of $\mathbf{R}$ is
\begin{small}
$$
  f_{\mathbf{R}}(\mathbf{R})=\frac{\displaystyle\prod_{i<j}^{N-1}(L_{i}+L_{j})}
  {|\boldsymbol{\Sigma}|^{K/2}}\sum_{t=0}^{\infty}\frac{h^{(t)}\left[\tr
  \left(\boldsymbol{\Sigma}^{-1}\mathbf{R}+\boldsymbol{\Omega}\right)\right]}{t!}
  \int_{V_{N-1,K}}\left[\tr\left(-2\boldsymbol{\mu}'\boldsymbol{\Sigma}^{-1}
  \mathbf{R}\mathbf{H}\right)\right]^{t}(\mathbf{H}d\mathbf{H}').
$$
\end{small}
The integral equals zero when $t$ is odd, then by \citet[eq.
(22)]{JAT64}
$$
  \int_{V_{N-1,K}}\left[\tr\left(-2\boldsymbol{\mu}'\boldsymbol{\Sigma}^{-1}
  \mathbf{R}\mathbf{H}\right)\right]^{2t}(\mathbf{H}d\mathbf{H}')=
  \frac{2^{N-1}\pi^{\frac{(N-1)K}{2}}}{\Gamma_{N-1}\left[\frac{K}{2}\right]}\sum_{\kappa}
  \frac{\left(\frac{1}{2}\right)_{t}4^{t}}{\left(\frac{1}{2}K\right)_{\kappa}}
  C_{\kappa}\left(\boldsymbol{\Omega}\boldsymbol{\Sigma}^{-1}\mathbf{R}^{2}\right).
$$
Noting that $\frac{\left(\frac{1}{2}\right)_{t}4^{t}}{(2t)!} = \frac{1}{t!}$,  so
$$
  f_{\mathbf{R}}(\mathbf{R})=\frac{\pi^{\frac{(N-1)K}{2}}\displaystyle
  \prod_{i<j}^{N-1}(L_{i}+L_{j})}{2^{-N+1}\Gamma_{N-1}\left[\frac{K}{2}\right]
  |\boldsymbol{\Sigma}|^{\frac{K}{2}}} \sum_{t=0}^{\infty}\sum_{\kappa}
  \frac{h^{(2t)}\left[\tr\left(\boldsymbol{\Sigma}^{-1}\mathbf{R}^{2} +
  \boldsymbol{\Omega}\right)\right]
  C_{\kappa}\left(\boldsymbol{\Omega}\boldsymbol{\Sigma}^{-1}\mathbf{R}^{2}\right)}
 {t!\left(\frac{1}{2}K\right)_{\kappa}}. \qed
$$

\section{Polar shape density}\label{sec:polarshape}

Now, observe that $\mathbf{R}:N-1\times N-1$, $\mathbf{R}>0$,
contains $(N-1)N/2$ different coordinates $(r_{ij}=r_{ji})$. Let
$v(\mathbf{R})$ the vector consisting of the different elements
$r_{ij}$, taken column by column.  Then the polar shape matrix
$\mathbf{W}$, can be written as:
$$
v(\mathbf{W})=\frac{1}{r}v(\mathbf{R}), \quad
r=\|\mathbf{R}\|=\sqrt{\tr \mathbf{R}^{2}}=\|\mathbf{Y}\|.
$$
Then by \citet{MR1982}, Theorem 2.1.3, p. 55:
$$
(dv(\mathbf{W}))=r^{m}\prod_{i=1}^{m}\sin
^{m-i}\boldsymbol{\Theta}_{i}\bigwedge_{i=1}^{m}d\boldsymbol{\Theta}_{i}\wedge
dr,\quad m=\frac{N(N-1)}{2}-1,
$$
which will denoted as
$$
(d\mathbf{W})=r^{m}J(\mathbf{u})\bigwedge_{i=1}^{m}d\boldsymbol{\Theta}_{i}\wedge
dr.
$$
Thus:
\begin{thm}\label{th:Polarshape}
The polar  shape density is
\begin{eqnarray*}
    f_{\mathbf{W}}(\mathbf{W})&=&\frac{2^{N-1}\pi^{\frac{(N-1)K}{2}}
    \displaystyle\prod_{i<j}^{N-1}(\lambda_{i}+\lambda_{j})J(\mathbf{u})}
    {\Gamma_{N-1}\left[\frac{K}{2}\right]|\boldsymbol{\Sigma}|^{\frac{K}{2}}}
    \sum_{t=0}^{\infty}\sum_{\kappa}\frac{C_{\kappa}\left(\boldsymbol{\Omega}
    \boldsymbol{\Sigma}^{-1}\mathbf{W}^{2}\right)}{t!\left(\frac{1}{2}K\right)_{\kappa}}
    \\&&\times\int_{0}^{\infty}r^{(N-1)^{2}+2t-1}h^{(2t)}\left[r^{2}\tr
    \boldsymbol{\Sigma}^{-1}\mathbf{W}^{2}+\tr\boldsymbol{\Omega}\right](dr).
\end{eqnarray*}
\end{thm}
\textit{Proof.}  The density of $\mathbf{R}$ is
\begin{small}
\begin{equation*}
   f_{\mathbf{R}}(\mathbf{R})=\frac{\pi^{\frac{(N-1)K}{2}}\displaystyle
  \prod_{i<j}^{N-1}(L_{i}+L_{j})} {2^{-N+1}\Gamma_{N-1}\left[\frac{K}{2}\right]
  |\boldsymbol{\Sigma}|^{\frac{K}{2}}} \sum_{t=0}^{\infty}\sum_{\kappa}
  \frac{h^{(2t)}\left[\tr\left(\boldsymbol{\Sigma}^{-1}\mathbf{R}^{2}
  +\boldsymbol{\Omega}\right)\right]
  C_{\kappa}\left(\boldsymbol{\Omega}\boldsymbol{\Sigma}^{-1}\mathbf{R}^{2}\right)}
  {t!\left(\frac{1}{2}K\right)_{\kappa}}.
\end{equation*}
\end{small}
Let be $\mathbf{W}(\mathbf{u})=\mathbf{R}/r$, then the joint density function of
$\mathbf{W}(\mathbf{u})$ and $r$ is given by
\begin{eqnarray*}\hspace{-1cm}
    f_{r,\mathbf{W}(\mathbf{u})}(r,\mathbf{W}(\mathbf{u}))&=&\frac{2^{N-1}\pi^{\frac{(N-1)K}{2}}
    \displaystyle\prod_{i<j}^{N-1}\left(r\left(\lambda_{i}+\lambda_{j}\right)\right)}
    {\Gamma_{N-1}\left[\frac{K}{2}\right]|\boldsymbol{\Sigma}|^{\frac{K}{2}}}\\
    && \times \ \sum_{t=0}^{\infty}\sum_{\kappa} \frac{h^{(2t)}
    \left[\tr\left(r^{2}\boldsymbol{\Sigma}^{-1}\mathbf{W}^{2} +
    \boldsymbol{\Omega}\right)\right]C_{\kappa}\left(r^{2}
    \boldsymbol{\Omega}\boldsymbol{\Sigma}^{-1}\mathbf{W}^{2}\right)}
    {t!\left(\frac{1}{2}K\right)_{\kappa}}r^{m}J(\mathbf{u}),
\end{eqnarray*}
with $m=N(N-1)/2-1$. Let $\lambda_{i}=\lambda_{i}(\mathbf{W})$ of $\mathbf{W}$, so if
$L_{i}=\lambda_{i}(\mathbf{R})$, thus $L_{i}=r\lambda_{i}$. Also note that:
\begin{enumerate}
    \item
    $C_{\kappa}\left(r^{2}\boldsymbol{\Omega}\boldsymbol{\Sigma}^{-1}\mathbf{W}^{2}\right) =
    r^{2t}C_{\kappa}\left(\boldsymbol{\Omega}\boldsymbol{\Sigma}^{-1}\mathbf{W}^{2}\right)$,
    \item
    $\displaystyle\prod_{i<j}^{N-1}r(\lambda_{i}+\lambda_{j})=r^{(N-1)(N-2)/2}\displaystyle
    \prod_{i<j}^{N-1}(\lambda_{i}+\lambda_{j})$,
    \item
    $h^{(2t)}\left[\tr\left(r^{2}\boldsymbol{\Sigma}^{-1}\mathbf{W}^{2} + \boldsymbol{\Omega}
    \right)\right]=h^{(2t)}\left[r^{2}\tr\boldsymbol{\Sigma}^{-1}\mathbf{W}^{2}+
    \tr\boldsymbol{\Omega}\right]$.
\end{enumerate}
Collecting powers of $r$ as $r^{m+2t+(N-1)(N-2)/2}=r^{(N-1)^{2}+2t-1}$ , the marginal
of $\mathbf{W}$ is
\begin{eqnarray*}
    f_{\mathbf{W}}(\mathbf{W})&=&\frac{2^{N-1}\pi^{\frac{(N-1)K}{2}}\displaystyle
    \prod_{i<j}^{N-1}\left(\lambda_{i}+\lambda_{j}\right)J(\mathbf{u})}
    {\Gamma_{N-1}\left[\frac{K}{2}\right]|\boldsymbol{\Sigma}|^{\frac{K}{2}}}
    \sum_{t=0}^{\infty}\sum_{\kappa}\frac{C_{\kappa}\left(\boldsymbol{\Omega}
    \boldsymbol{\Sigma}^{-1}\mathbf{W}^{2}\right)}{t!\left(\frac{1}{2}K\right)_{\kappa}}\\
    && \times \ \int_{0}^{\infty}r^{(N-1)^{2}+2t-1}h^{(2t)}
    \left[r^{2}\tr\boldsymbol{\Sigma}^{-1}\mathbf{W}^{2}+\tr\boldsymbol{\Omega}\right](dr).\qed
\end{eqnarray*}

\begin{rem}\label{rem:Polarnoreflection}
Given that $\mathbf{H}\in V_{N-1,K}$, we cannot classify the polar shape densities by
including or excluding reflections as in the QR shape distribution cases.
\end{rem}

\section{Central case}\label{sec:central}

The central case of the elliptical polar shape densities follows
easily:

\begin{cor}\label{cor:Polarcentralsizeandshape}
The central polar size-and-shape density is given by
\begin{equation*}
    f_{\mathbf{R}}(\mathbf{R})=\frac{2^{N-1}\pi^{\frac{(N-1)K}{2}}\displaystyle
    \prod_{i<j}^{N-1}(L_{i}-L_{j})}{\Gamma_{n}\left[\frac{K}{2}\right]
    |\boldsymbol{\Sigma}|^{\frac{K}{2}}} h\left[\tr \boldsymbol{\Sigma}^{-1}
    \mathbf{R}^{2}\right]
\end{equation*}
\end{cor}
\textit{Proof.} Just take $\boldsymbol{\mu}=0$ in Theorem \ref{th:Polarsizeandshape}
and use $h^{(0)}(\cdot)=h(\cdot)$.\qed

And finally, we have that

\begin{cor}\label{cor:Polarcentralinvariance}
The central polar shape density is invariant under the elliptical family and it is
given by
\begin{equation}
    f_{\mathbf{W}}(\mathbf{W})=\frac{2^{N-2}\pi^{\frac{(N-1)(K-N)}{2}}
    \Gamma\left[(N-1)^{2}\right]}{\Gamma_{N-1}\left[\frac{K}{2}\right]
    |\boldsymbol{\Sigma}|^{\frac{K}{2}}}\prod_{i<j}(\lambda_{i}+\lambda_{j})J(\mathbf{u})
    \left(\tr\boldsymbol{\Sigma}^{-1}\mathbf{W}^{2}\right)^{-\frac{(N-1)^{2}}{2}}
\end{equation}
\end{cor}
\textit{Proof.} It is straightforward from Theorem
\ref{th:Polarshape}. Take $\boldsymbol{\mu}=0$, and use
$h^{(0)}(\cdot)=h(\cdot)$, then
\begin{equation*}
    f_{\mathbf{W}}(\mathbf{W})=\frac{2^{N-1}\pi^{\frac{(N-1)K}{2}}\displaystyle
    \prod_{i<j}^{N-1}(\lambda_{i}+\lambda_{j})}{\Gamma_{N-1}\left[\frac{K}{2}\right]
    |\boldsymbol{\Sigma}|^{\frac{K}{2}}}J(\mathbf{u}) \int_{0}^{\infty}r^{(N-1)^{2}-1}
    h\left[r^{2}\tr\boldsymbol{\Sigma}^{-1}\mathbf{W}^{2}\right](dr)
\end{equation*}
Let be
$s=\left(\tr\boldsymbol{\Sigma}^{-1}\mathbf{W}^{2}\right)^{\frac{1}{2}}r$,
so
$ds=\left(\tr\boldsymbol{\Sigma}^{-1}\mathbf{W}^{2}\right)^{\frac{1}{2}}(dr)$,
and
\begin{small}
\begin{eqnarray*}
    &&\int_{0}^{\infty}\left(\frac{s}{\left(\tr\boldsymbol{\Sigma}^{-1}\mathbf{W}^{2}
    \right)^{\frac{1}{2}}}\right)^{(N-1)^{2}-1} h\left(s^{2}\right)\frac{ds}{\left(\tr
    \boldsymbol{\Sigma}^{-1}\mathbf{W}^{2}\right)^{\frac{1}{2}}}\hspace{4cm}\\
    && \hspace{4cm} =\left(\tr\boldsymbol{\Sigma}^{-1} \mathbf{W}^{2}\right)^{-\frac{(N-1)^{2}}{2}}
    \int_{0}^{\infty}s^{N(N-1)-1}h\left(s^{2}\right)(ds)\\
    && \hspace{4cm}  =\frac{\Gamma\left[(N-1)^{2}\right]}{2\pi^{\frac{(N-1)^{2}}{2}}
    \left(\tr\boldsymbol{\Sigma}^{-1}\mathbf{W}^{2}\right)^{\frac{(N-1)^{2}}{2}}}.
\end{eqnarray*}
\end{small}
Then
\begin{equation*}
    f_{\mathbf{W}}(\mathbf{W})=\frac{\pi^{(N-1)(K-N)/2}\Gamma\left[(N-1)^{2}\right]}{2^{-N+2}
    \Gamma_{N-1}\left(\frac{K}{2}\right)|\boldsymbol{\Sigma}|^{\frac{K}{2}}}
    \displaystyle\prod_{i<j}^{N-1}\left(\lambda_{i}+\lambda_{j}\right)
    J(\mathbf{u})\left(\tr\boldsymbol{\Sigma}^{-1}\mathbf{W}^{2}\right)^{-\frac{(N-1)^{2}}{2}}.
    \qed
\end{equation*}

\section{Some particular models}\label{sec:particularmodelsexample}

Finally, we give explicit shapes densities for some elliptical models.

The Kotz type I model is given by
$$
  h(y)=\frac{R^{T-1+\frac{K(N-1)}{2}}\Gamma\left(\frac{K(N-1)}{2}\right)}{\pi^{K(N-1)/2}
  \Gamma\left(T-1+\frac{K(N-1)}{2}\right)}y^{T-1}\exp\{-Ry\}.
$$
So, the corresponding $k$-th derivative follows from
\begin{eqnarray*}
  &&\frac{d^{k}}{dy^{k}}y^{T-1}\exp\{-Ry\}= \\
  && \qquad (-R)^{k}y^{T-1}\exp\{-Ry\}\left\{1+\sum_{m=1}^{k}\binom{k}{m}
  \left[\prod_{i=0}^{m-1}(T-1-i)\right](-Ry)^{-m}\right\},
\end{eqnarray*}
see \citet{Caro2009}.

It is of interest the  normal case, i.e. when $T=1$ and $R=\frac{1}{2}$, here the
derivation is straightforward from the general density.

The required derivative follows easily, it is,
$$
  h^{(k)}(y)=\frac{R^{\frac{K(N-1)}{2}}}{\pi^{\frac{K(N-1)}{2}}}(-R)^{k}\exp\{-Ry\}
$$
and replacing
\begin{eqnarray*}
    &&\int_{0}^{\infty}r^{(N-1)^{2}+2t-1}h^{(2t)}\left[r^{2}\tr\boldsymbol{\Sigma}^{-1}
    \mathbf{W}^{2}+\tr\boldsymbol{\Omega}\right]dr\\
    && \hspace{2cm}=\frac{R^{\frac{K(N-1)}{2}+2t}
    \etr\{-R\boldsymbol{\Omega}\}}{2\pi^{\frac{K(N-1)}{2}}\left(\tr
    R\boldsymbol{\Sigma}^{-1}\mathbf{W}^{2}\right)^{\frac{(N-1)^{2}}{2}+t}}
    \Gamma\left[\frac{(N-1)^{2}}{2}+t\right],
\end{eqnarray*}
in

\begin{eqnarray*}
    f_{\mathbf{W}}(\mathbf{W})&=&\frac{2^{N-1}\pi^{\frac{(N-1)K}{2}}\displaystyle
    \prod_{i<j}^{N-1}\left(\lambda_{i}+\lambda_{j}\right)J(\mathbf{u})}
    {\Gamma_{N-1}\left[\frac{K}{2}\right]|\boldsymbol{\Sigma}|^{\frac{K}{2}}}
    \sum_{t=0}^{\infty}\sum_{\kappa}\frac{C_{\kappa}\left(\boldsymbol{\Omega}
    \boldsymbol{\Sigma}^{-1}\mathbf{W}^{2}\right)}{t!\left(\frac{1}{2}K\right)_{\kappa}}\\
    &&\times\int_{0}^{\infty}r^{(N-1)^{2}+2t-1}h^{(2t)}\left[r^{2}\tr
    \boldsymbol{\Sigma}^{-1}\mathbf{W}^{2}+\tr\boldsymbol{\Omega}\right]dr,
\end{eqnarray*}
we have proved that
\begin{cor}\label{cor:polarshapeNORMAL}
The Gaussian polar  shape density is
\begin{eqnarray*}
    f_{\mathbf{W}}(\mathbf{W}) &=&\frac{2^{N-2}J(\mathbf{u})\displaystyle
    \prod_{i<j}^{N-1}\left(\lambda_{i}+\lambda_{j}\right)} {R^{-\frac{K(N-1)}{2}}
    \Gamma_{N-1}\left[\frac{K}{2}\right]|\boldsymbol{\Sigma}|^{\frac{K}{2}}}
    \frac{\etr\{-R\boldsymbol{\Omega}\}}{\left(\tr R \boldsymbol{\Sigma}^{-1}
    \mathbf{W}^{2}\right)^{\frac{(N-1)^{2}}{2}}}\\
    && \sum_{t=0}^{\infty}\frac{\Gamma\left[\frac{(N-1)^{2}}{2}+t\right]}{t! \left(\tr R
    \boldsymbol{\Sigma}^{-1}\mathbf{W}^{2}\right)^{t}}\sum_{\kappa}\frac{C_{\kappa}\left(R^{2}
    \boldsymbol{\Omega}\boldsymbol{\Sigma}^{-1}\mathbf{W}^{2}\right)}{\left(\frac{1}{2}
    K\right)_{\kappa}}.
\end{eqnarray*}
\end{cor}

Finally, we propose the result for the Kotz type I model
$$
  h(y)=\frac{R^{T-1+\frac{K(N-1)}{2}}\Gamma\left(\frac{K(N-1)}{2}\right)}{\pi^{K(N-1)/2}
  \Gamma\left(T-1+\frac{K(N-1)}{2}\right)}y^{T-1}\exp\{-Ry\}.
$$

\begin{cor}\label{cor:polarshapeKotz}
The Kotz type I polar  shape density is
\begin{eqnarray*}
    f_{\mathbf{W}}(\mathbf{W})&=&\frac{2^{N-1}\displaystyle\prod_{i<j}^{N-1}
    \left(\lambda_{i}+\lambda_{j}\right)J(\mathbf{u})} {\Gamma_{N-1}
    \left[\frac{K}{2}\right]| \boldsymbol{\Sigma}|^{\frac{K}{2}}\etr
    \{R\boldsymbol{\Omega}\}} \sum_{t=0}^{\infty}\sum_{\kappa}
    \frac{C_{\kappa}\left(\boldsymbol{\Omega}\boldsymbol{\Sigma}^{-1}
    \mathbf{W}^{2}\right)}{t!\left(\frac{1}{2}K\right)_{\kappa}}\\
    && \times \ \frac{R^{T-1+\frac{K(N-1)}{2}+2t}\Gamma\left(\frac{K(N-1)}{2}\right)}
    {\Gamma\left(T-1+\frac{K(N-1)}{2}\right)}\\
    && \times \ \left\{\sum_{i=0}^{\infty}\frac{1}{i!}\prod_{u=0}^{i-1}(T-1-u)
    \times\frac{(\tr\boldsymbol{\Omega})^{T-1-i}\Gamma\left[\frac{(N-1)^{2}}{2}+i+t\right]}
    {2R^{\frac{(N-1)^{2}}{2}+i+t}(\tr \boldsymbol{\Sigma}^{-1}
    \mathbf{W}^{2})^{\frac{(N-1)^{2}}{2}+t}}\right.\\
    && + \ \sum_{m=1}^{2t}\binom{2t}{m}\left[\prod_{i=0}^{m-1}(T-1-i)\right](-R)^{-m}\\
    && \times \  \left.\sum_{i=0}^{\infty}\frac{1}{i!}\prod_{u=0}^{i-1}(T-1-m-u)
    \frac{(\tr\boldsymbol{\Omega})^{T-1-m-i} \Gamma\left[\frac{(N-1)^{2}}{2}+t+i\right]}
    {2R^{\frac{(N-1)^{2}}{2}+i+t}(\tr\boldsymbol{\Sigma}^{-1}
    \mathbf{W}^{2})^{\frac{(N-1)^{2}}{2}+t}}\right\}.
\end{eqnarray*}
\end{cor}

\textit{Proof.} The corresponding $k$-th derivative follows from
\begin{eqnarray*}
  \frac{d^{k}}{dy^{k}}y^{T-1}\exp\{-Ry\}= \hspace{8.5cm}\\
  (-R)^{k}y^{T-1}\exp\{-Ry\}\left\{1+\sum_{m=1}^{k}\binom{k}{m}
  \left[\prod_{i=0}^{m-1}(T-1-i)\right](-Ry)^{-m}\right\},
\end{eqnarray*}
(see \citet{Caro2009}), and the corresponding polar  shape density
is obtained after some simplification like
\begin{eqnarray*}\hspace{-1cm}
    f_{\mathbf{W}}(\mathbf{W})&=&\frac{2^{N-1}\pi^{\frac{(N-1)K}{2}}\displaystyle\prod_{i<j}^{N-1}
    \left(\lambda_{i}+\lambda_{j}\right)J(\mathbf{u})} {\Gamma_{N-1}
    \left[\frac{K}{2}\right]| \boldsymbol{\Sigma}|^{\frac{K}{2}}} \sum_{t=0}^{\infty}
    \sum_{\kappa}\frac{C_{\kappa}\left(\boldsymbol{\Omega}
    \boldsymbol{\Sigma}^{-1}\mathbf{W}^{2}\right)}{t!\left(\frac{1}{2}K\right)_{\kappa}}\\
    && \times\int_{0}^{\infty}r^{(N-1)^{2}+2t-1}h^{(2t)}\left[r^{2} \tr
    \boldsymbol{\Sigma}^{-1} \mathbf{W}^{2} + \tr\boldsymbol{\Omega}\right]dr\\
    &=& \frac{2^{N-1}\displaystyle\prod_{i<j}^{N-1}\left(\lambda_{i}+\lambda_{j}
    \right)J(\mathbf{u})} {\Gamma_{N-1}\left[\frac{K}{2}\right]
    |\boldsymbol{\Sigma}|^{\frac{K}{2}}\etr\{R\boldsymbol{\Omega}\}}
    \sum_{t=0}^{\infty}\sum_{\kappa}\frac{C_{\kappa}
    \left(\boldsymbol{\Omega}\boldsymbol{\Sigma}^{-1}
    \mathbf{W}^{2}\right)}{t!\left(\frac{1}{2}K\right)_{\kappa}}\\
    && \times \ \frac{R^{T-1+\frac{K(N-1)}{2}+2t}\Gamma\left(\frac{K(N-1)}{2}\right)}
    {\Gamma\left(T-1+\frac{K(N-1)}{2}\right)} \\
    &&  \times \ \left\{\sum_{i=0}^{\infty}\frac{1}{i!}\prod_{u=0}^{i-1}(T-1-u) \times
    \frac{(\tr\boldsymbol{\Omega})^{T-1-i}\Gamma\left[\frac{(N-1)^{2}}{2}+i+t\right]}
    {2R^{\frac{(N-1)^{2}}{2}+i+t}(\tr \boldsymbol{\Sigma}^{-1}
    \mathbf{W}^{2})^{\frac{(N-1)^{2}}{2}+t}}\right.\\
    && + \ \sum_{m=1}^{2t}\binom{2t}{m}\left[\prod_{i=0}^{m-1}(T-1-i)\right](-R)^{-m}\\
    &&  \times \ \left.\sum_{i=0}^{\infty}\frac{1}{i!} \prod_{u=0}^{i-1}(T-1-m-u)
    \frac{(\tr\boldsymbol{\Omega})^{T-1-m-i} \Gamma\left[\frac{(N-1)^{2}}{2}+t+i\right]}
    {2R^{\frac{(N-1)^{2}}{2}+i+t}(\tr\boldsymbol{\Sigma}^{-1}\mathbf{W}^{2})^{
    \frac{(N-1)^{2}}{2}+t}}\right\}\qed
\end{eqnarray*}

\subsection{Example: Mouse Vertebra}\label{sub:mouse}
This experiment has been  studied in the Gaussian case  by \citet{DM98} and some
references there in. Here we study this data under three different models, the usual
normal and two non normal Kotz models, the best distribution will determined by
applying a modified BIC criterion.

We start with  the isotropic version  of the Gaussian polar density given in
corollary \ref{cor:polarshapeNORMAL}
\begin{eqnarray*}
    f_{\mathbf{W}}(\mathbf{W}) &=&\frac{J(\mathbf{u})\displaystyle\prod_{i<j}^{N-1}
    \left(\lambda_{i}+\lambda_{j}\right)\etr\left(-\frac{\boldsymbol{\mu}'\boldsymbol{\mu}}
    {2\boldsymbol{\Sigma}^{2}}\right)} {2^{\frac{K(N-1)}{2}-\frac{(N-1)^{2}}{2}-N+2}
    \Gamma_{N-1}\left[\frac{K}{2}\right]\boldsymbol{\Sigma}^{K(N-1)-(N-1)^{2}}}\\
    && \qquad \times \ \sum_{t=0}^{\infty}\frac{\Gamma\left[\frac{(N-1)^{2}}{2}+t\right] }{t!}
    \sum_{\kappa}\frac{C_{\kappa}\left(\frac{1}{2\boldsymbol{\Sigma}^{2}}
    \boldsymbol{\mu}'\mathbf{W}^{2}\boldsymbol{\mu}\right)}{\left(\frac{1}{2}K\right)_{\kappa}}.
\end{eqnarray*}

Consider now two Kotz models from corollary
\ref{cor:polarshapeKotz}:  for $T=2$ and $R=\frac{1}{2}$, we obtain
after some simplification that

\begin{eqnarray*}
    f_{\mathbf{W}}(\mathbf{W})&=&\frac{2^{\frac{N^{2}-K(N-1)-1}{2}}J(\mathbf{u})\displaystyle
    \prod_{i<j}^{N-1}\left(\lambda_{i}+\lambda_{j}\right)}
    {\boldsymbol{\Sigma}^{-(N-1)^{2}+K(N-1)}K(N-1) \Gamma_{N-1}\left[\frac{K}{2}\right]}
    \etr\left(-\frac{\boldsymbol{\mu}'\boldsymbol{\mu}}{2\boldsymbol{\Sigma}^{2}}\right)\\
    && \times\sum_{t=0}^{\infty}\frac{\left(-2t+\tr\left(\frac{\boldsymbol{\mu}'\boldsymbol{\mu}}
    {2\boldsymbol{\Sigma}^{2}}\right)\right)\Gamma\left[\frac{(N-1)^{2}}{2}+t\right]
    +\Gamma\left[\frac{(N-1)^{2}}{2}+t+1\right]}{t!}\\
    &&\times\sum_{\kappa}\frac{C_{\kappa}\left(\frac{1}{2\boldsymbol{\Sigma}^{2}}
    \boldsymbol{\mu}'\mathbf{W}^{2}\boldsymbol{\mu}\right)}{\left(\frac{1}{2}K\right)_{\kappa}};
\end{eqnarray*}

and, for $T=3$ and $R=\frac{1}{2}$, the corresponding non Gaussian
isotropic model is given by

\begin{eqnarray*}
    &&f_{\mathbf{W}}(\mathbf{W})=\frac{2^{\frac{N^{2}-K(N-1)+1}{2}}
    \boldsymbol{\Sigma}^{(N-1)^{2}-K(N-1)}J(\mathbf{u})\displaystyle
    \prod_{i<j}^{N-1}\left(\lambda_{i}+\lambda_{j}\right)}
    {K(N-1)(K(N-1)+2)\Gamma_{N-1}\left[\frac{K}{2}\right]}\etr
    \left(-\frac{\boldsymbol{\mu}'\boldsymbol{\mu}}{2\boldsymbol{\Sigma}^{2}}\right)\\
    && \times\sum_{t=0}^{\infty}\frac{1}{t!}\left\{\left(-2t+4t^{2}-4t \tr
    \left(\frac{\boldsymbol{\mu}'\boldsymbol{\mu}}{2\boldsymbol{\Sigma}^{2}}\right)
    +\tr^{2}\left(\frac{\boldsymbol{\mu}'\boldsymbol{\mu}}{2\boldsymbol{\Sigma}^{2}}
    \right)\right)\Gamma\left[\frac{(N-1)^{2}}{2}+t\right]\right.\\
    &&\left.+\left(-4t+2\tr\left(\frac{\boldsymbol{\mu}'
    \boldsymbol{\mu}}{2\boldsymbol{\Sigma}^{2}}\right)\right)
    \Gamma\left[\frac{(N-1)^{2}}{2}+t+1\right]+
    \Gamma\left[\frac{(N-1)^{2}}{2}+t+2\right]\right\}\\
    &&\times\sum_{\kappa}\frac{C_{\kappa}\left(\frac{1}{2\boldsymbol{\Sigma}^{2}}
    \boldsymbol{\mu}'\mathbf{W}^{2}\boldsymbol{\mu}\right)}{\left(\frac{1}{2}K\right)_{\kappa}}.
\end{eqnarray*}
We contrast these three models via the modified BIC criterion, they will be applied
to the data of two groups (small and large) of mouse vertebra, an experiment very
detailed in \citet{DM98}.

The likelihood based on the exact densities require the computation of the above
series, a carefully comparison with the known hypergeometric of one matrix argument
indicates that these distributions can be obtained by a suitable modification of the
algorithms of \citet{KE06}.

In order to decide which the elliptical model is the best one, different criteria
have been employed for the model selection. We shall consider a modification of the
BIC statistic as discussed in \citet{YY07}, and which was first achieved by
\citet{ri:78} in a coding theory framework. The modified BIC is given by:
$$
  BIC^{*}=-2\mathfrak{L}(\widetilde{\boldsymbol{\boldsymbol{\boldsymbol{\mu}}}},
  \widetilde{\boldsymbol{\Sigma}}^{2},h) +n_{p}(\log(n+2)-\log 24),
$$
where $\mathfrak{L}(\widetilde{\boldsymbol{\boldsymbol{\boldsymbol{\mu}}}},
\widetilde{\boldsymbol{\Sigma}}^{2},h)$ is the maximum of the log-likelihood
function, $n$ is the sample size and $n_{p}$ is the number of parameters to be
estimated for each particular shape density.

As proposed by \citet{kr:95} and \citet{r:95}, the following selection criteria have
been employed for the model selection.

\begin{table}[ht]  \centering \caption{Grades of evidence
corresponding to values of the $BIC^{*}$ difference.}\label{table2}
\medskip
\renewcommand{\arraystretch}{1}
\begin{center}
  \begin{tabular}{cl}
    \hline
    $BIC^{*}$ difference & Evidence\\
    \hline
    0--2 & Weak\\
    2--6 & Positive \\
    6--10 & Strong\\
    $>$ 10 & Very strong\\
    \hline
  \end{tabular}
\end{center}
\end{table}

In order to apply the above densities we need to restrict the number of landmarks in
such way that $\min(K,N-1)=N-1$, so in the mouse vertebra data we must select 3
landmarks of the original 6 points. In the following example we consider the
landmarks 1, 2 and 6 which corresponds to the widest part of the vertebra.

The maximum likelihood estimators  for location  and scale parameters  associated
with the small and large groups are summarised in the following table:

\begin{table}[ht]  \centering \caption{Maximum likelihood estimators.}\label{table3}
\medskip
\renewcommand{\arraystretch}{1}
\begin{center}
\begin{small}\label{Tab4MR}
\begin{tabular}{|c|c|c|c|c|c|c|}
  \hline\hline
  Group& $BIC^{*}$ & $\widetilde{\mu}_{11}$ &$\widetilde{\mu}_{12}$& $\widetilde{\mu}_{21}$ & $\widetilde{\mu}_{22}$
   & $\widetilde{\sigma}^{2}$   \\
   & $\build{K:T=2}{G}{K:T=3}$&&&&&\\
  \hline\hline
  Small & $\build{-169.0557}{-129.9719}{-139.2333}$ & $\build{-6.3553}{-4.5092}{-5.6873}$ &
$\build{-54.5119}{-52.4558}{-50.9245}$
   & $\build{-18.9629}{-18.2941}{-20.4403}$
  &$\build{3.8065}{3.0990}{4.8584}$    &$\build{28.1767}{32.0686}{27.0289}$ \\
  \hline
\end{tabular}

\medskip

\begin{tabular}{|c|c|c|c|c|c|c|}
   \hline\hline
   Group& $BIC^{*}$ & $\widetilde{\mu}_{11}$ &$\widetilde{\mu}_{12}$& $\widetilde{\mu}_{21}$
   & $\widetilde{\mu}_{22}$ & $\widetilde{\sigma}^{2}$ \\
   & $\build{K:T=3}{G}{K:T=2}$&&&&&\\
   \hline\hline
   Large& $\build{-204.2375}{-128.2599}{-167.2111}$ &$\build{-20.1702}{-28.0776}{-26.0985}$
   &$\build{-95.3134}{-73.1440}{-82.8228}$ & $\build{-41.9462}{-32.0698}{-36.3675}$
   & $\build{9.8279 }{13.1131}{12.3299}$   & $\build{84.3441}{75.7233}{75.0727}$  \\
  \hline
\end{tabular}
\end{small}
\end{center}
\end{table}

According to the modified BIC criterion, the Kotz model with parameters $T=2$,
$R=\frac{1}{2}$ and $s=1$ is the most appropriate for the small group, instead the
Kotz distribution with parameters $T=3$, $R=\frac{1}{2}$ and $s=1$ models the large
group. There is a very strong difference between these models and the classical
normal in this experiment.

Let $\boldsymbol{\boldsymbol{\boldsymbol{\mu}}}_{1}$ and
$\boldsymbol{\boldsymbol{\boldsymbol{\mu}}}_{2}$ be the mean polar shape of the small
and large groups, respectively. We test equal mean shape under the best models, and
the likelihood ratio (based on $-2\log\Lambda\approx\chi_{4}^{2}$) for the test
$H_{0}:\boldsymbol{\boldsymbol{\boldsymbol{\mu}}}_{1} =
\boldsymbol{\boldsymbol{\boldsymbol{\mu}}}_{2}$ vs
$H_{a}:\boldsymbol{\boldsymbol{\boldsymbol{\mu}}}_{1} \neq\
\boldsymbol{\boldsymbol{\boldsymbol{\mu}}}_{2}$, provides the p-value $0.99$, which
means that there extremely evidence that the mean shapes of the two groups are equal.

For any elliptical model we can obtain the polar shape density, however a nontrivial
problem appears, the $2t$-th derivative of the generator model, which can be seen as
a partition theory problem. For The general case of a Kotz model ($s\neq1$), and
another models like Pearson II and VII, Bessel, Jensen-logistic, we can use formulae
for these derivatives given by \citet{Caro2009}. The resulting densities have again a
form of a generalised series of zonal polynomials which can be computed efficiently
after some modification of existing works for hypergeometric series (see
\citet{KE06}), thus the inference over an exact density can be performed, avoiding
the use of any asymptotic distribution, and the initial transformation avoids the
invariant polynomials of \citet{D80}, which at present seems can not be computable
for large degrees.

\section*{Acknowledgment}

This research work was supported  by University of Medellin (Medellin, Colombia) and
Universidad Aut\'onoma Agraria Antonio Narro (M\'{e}xico),  joint grant No. 469,
SUMMA group. Also, the first author was partially supported  by CONACYT - M\'exico,
research grant no. \ 138713 and IDI-Spain, Grants No. \ FQM2006-2271 and
MTM2008-05785.


\begin{thebibliography}{}

\bibitem[Cadet(1996)]{c:96}
    Cadet, A. (1996)
    Polar coordinates in $\Re^{np}$; Application to the computation of Wishart and Beta
    laws.
    \emph{Sankhy$\bar{a}$ A,} \textbf{58}, 101-114.

\bibitem[Caro-Lopera \textit{et al} (2009)]{Caro2009}
    Caro-Lopera, F. J., D\'{\i}az-Garc\'{\i}a J. A. and Gonz\'{a}lez-Far\'{\i}as, G. (2009)
    Noncentral elliptical configuration density.
    \emph{J. Multivariate Anal.}, \textbf{101}(1), Pages 32-43.

\bibitem[Davis (1980)]{D80}
    Davis, A. W. (1980).
    Invariant polynomials with two matrix arguments, extending the zonal polynomials.
    In \textit{Multivariate analysis, V} (P.R. Krishnaiah, Ed.), 287--299, North-Holland,
    Amsterdam-New York.

\bibitem[D\'{\i}az-Garc\'{\i}a \emph{et al.}(1997)]{dgm:97}
    D\'{\i}az-Garc\'{\i}a, J. A., Guti\'errez- J\'aimez, R. and Mardia, K. V. (1997)
    Wishart and Pseudo-Wishart distributions and some applications to shape theory.
    \emph{J. Multivariate Anal.}, \textbf{63}, 73-87.

\bibitem[D\'{\i}az-Garc\'{\i}a \textit{et al}.(2003)]{dgr:03}
    D\'{\i}az-Garc\'{\i}a, J. A., Guti\'errez-J\'aimez, R. and Ramos, R. (2003)
    Size-and-Shape Cone, Shape Disk and Configuration Densities for the Elliptical Models.
    \emph{Braz. J. Probab. Stat.}, \textbf{17}, 135-146.

\bibitem[Dryden and Mardia (1998)]{DM98}
    Dryden, I. L. and Mardia, K. V. (1998)
    \emph{Statistical shape analysis}.
    John Wiley and Sons, Chichester, 1998.

\bibitem[Fang and Zhang (1990)]{fz:90}
    Fang, K. T.  and Zhang. Y. T.(1990)
    \emph{Generalized  Multivariate Analysis}.
    Science Press, Springer-Verlag, Beijing.

\bibitem[Goodall and Mardia (1993)]{GM93}
    Goodall, C. R. and Mardia. K. V. (1993)
    Multivariate Aspects of Shape Theory.
    \textit{Ann. Statist.}, \textbf{21}, 848--866.

\bibitem[Goodall(1991)]{g:91}
    Goodall, C. R. (1991).
    Procustes methods in the statistical analysis of shape (with discussion).
    \textit{J. Roy. Statist. Soc. Ser. B.}, \textbf{53}, 285-339.


\bibitem[Gupta and Varga(1993)]{gv:93}
    Gupta, A. K. and Varga, T. (1993)
    \textit{Elliptically Contoured Models in Statistics}.
    Kluwer Academic Publishers, Dordrecht.

\bibitem[James(1964)]{JAT64}
    James A. T. (1964).
    Distributions of matrix variate and latent roots derived from normal samples.
    \emph{Ann. Math. Statist.}, \textbf{35}, 475--501.

\bibitem[Kass and Raftery (1995)]{kr:95}
    R. E. Kass, R. E. and Raftery, A. E. (1995)
    Bayes factor.
    \emph{J. Amer. Statist. Soc.}, \textbf{90}, 773--795.


\bibitem[Koev and Edelman (2006)]{KE06}
    Koev, P. and Edelman A. (2006)
    The efficient evaluation of the hypergeometric function of a matrix argument.
    \textit{Math. Comp.} \textbf{75}, 833-846.

\bibitem[Le and Kendall (1993)]{LK93}
    Le, H. L. and Kendall, D. G. (1993)
    The Riemannian structure of Euclidean spaces: a novel environment for statistics.
    \emph{Ann. Statist.}, \textbf{21}, 1225--1271.

\bibitem[Muirhead(1982)]{MR1982}
    Muirhead R. J. (1982)
    \textit{Aspects of multivariate statistical theory}.
    Wiley Series in Probability and Mathematical Statistics.
    John Wiley \& Sons, Inc.,  New York.

\bibitem[Raftery(1995)]{r:95}
    Raftery, A. E. (1995)
    Bayesian model selection in social research.
    \emph{Sociological Methodology}, \textbf{25}, 111--163.

\bibitem[Rissanen(1978)]{ri:78}
    Rissanen, J. (1978)
    Modelling by shortest data description.
    \emph{Automatica}, \textbf{14}, 465--471.

\bibitem[Yang and Yang(2007)]{YY07}
    Yang, Ch. Ch. and Yang, Ch. Ch. (2007)
    Separating latent classes by information criteria.
    \emph{J. Classification}, \textbf{24}, 183--203.
\end{thebibliography}
\end{document}